\documentclass[11pt]{article}
\usepackage{eurosym}
\usepackage{graphicx}
\usepackage{enumerate}
\usepackage{amsmath}
\usepackage{amsfonts}
\usepackage{amssymb}
\usepackage{amsfonts}
\usepackage{amssymb,amsmath,amsthm,epsfig,graphicx,natbib,subfigure,hyperref}
\usepackage[margin=1.3in]{geometry}
\usepackage{caption2}
\usepackage{color}
\usepackage{enumerate}
\usepackage{booktabs}
\usepackage{supertabular}
\usepackage{float}

\def\blackbox{\hfill {\vrule height3pt width4pt depth2pt}}
\def\Box{\hfill \framebox(5.25,5.25){}}

\bibpunct[: ]{(}{)}{;}{a}{}{,}
\newtheorem{thm}{Theorem}
\newtheorem{pro}{Proposition}
\newtheorem{lem}{Lemma}

\newtheorem{corr}{Corollary}
\newtheorem{defin}{Definition}
\newtheorem{ob}{Observation}

\newenvironment{thm-prf}{\begin{thm} \nopagebreak}{\end{thm}}
\newenvironment{pro-prf}{\begin{pro} \nopagebreak}{\end{pro}}

\newenvironment{lem-prf}{\begin{lem} \nopagebreak}{\end{lem}}

\newenvironment{corr-prf}{\begin{corr} \nopagebreak}{\end{corr}}

\newcommand {\be}{\begin{equation}}
\newcommand {\ee}{\end{equation}}
\begin{document}

\title{Information Management for Decentralized Energy Storages under Market
Uncertainties}
\date{\today }
\author{Qiao-Chu~He, Yun~Yang, and~Baosen~Zhang \thanks{%
Q. He is with the Department of Systems Engineering and Engineering
Management, University of North Carolina, Charlotte, NC, 28223 USA. Email:
qhe4@uncc.edu.} \thanks{%
Y. Yang is with the Statistics Department, Florida State University, Tallahassee, FL, 32306 USA. Email:
yyang@stat.fsu.edu.} \thanks{%
B. Zhang is with the Electrical Engineering Department, University of Washington, Seattle, WA 98195 USA. Email:
zhangbao@uw.edu.}}

\maketitle

\begin{abstract}
\linespread{0.9} In this paper, we propose a model of decentralized energy
storages, who serve as instruments to shift energy supply intertemporally.
From storages' perspective, we investigate their optimal buying or selling
decisions under market uncertainty. The goal of this paper is to understand
the economic value of future market information, as energy storages mitigate
market uncertainty by forward-looking strategies.

At a system level, we evaluate different information management policies to
coordinate storages' actions and improve their aggregate profitability: (1)
providing a publicly available market forecasting channel; (2) encouraging
decentralized storages to share their private forecasts with each other; (3)
releasing additional market information to a targeted subset of storages
exclusively. We highlight the perils of too much market information
provision and advice on exclusiveness of market forecast channel.

\end{abstract}

%\author{Qiao-Chu He\thanks{%
%University of California at Berkeley, 4176 Etcheverry Hall, Berkeley, CA;
%e-mail:\ heqc0425@berkeley.edu.}}

%\author{}

%\and Ying-Ju Chen\thanks{%
%Hong Kong University of Science and Technology, LSK 4035, Clear Water Bay, Kowloon, Hong Kong; e-mail:\ imchen@ust.hk.}}

\setlength{\textwidth}{15cm} \setlength{\textheight}{21.5cm}

\linespread{1.0}

\psfull

\section{Introduction}

Energy storages serve as instruments for \emph{energy supply shift} by
storing excess renewable energy for the future. As the power industry
transitions from a regulated towards a more competitive market environment,
storage devices have incentives to be charged when prices are low
(corresponding to low residual demand), and discharges when prices are
higher (peak demand). The economic feasibility of such delicate operations
requires optimal utilization of decentralized storage devices, wherein
storage devices pursue their own objective (e.g., maximize profits or
minimize costs) \citep{Sarker2015}. Tesla Powerwall is among the most
renowned examples of such applications \citep{Tesla2016}.

The application scenarios we focus on in this paper are when energy storages
are integrated with renewable energy (e.g., wind and solar power)
generation. For example, energy supply shift is necessary when the peak in
solar power supply is during the daytime whereas the peak demand is during
the night. In this case, energy storages are placed near the wind and solar
power production sites to smooth generation output before connecting to
aggregator and feeding to the grid \citep{Grothoff2015}. However, a
fundamental problem in this integration is so called ``merit order effect":
The supply of renewable energy has negligible marginal costs and in turn
reduces the spot equilibrium price \citep{acemoglu2015competition}. For
example, since the price of power is expected to be lower during periods
with high wind than in periods with low wind, the intermittency of renewable
energy generation leads to market variability and uncertainty, i.e., energy
prices fluctuates dynamically over time.

Facing these challenges, energy storages have to employ market forecasts to
improve buying and selling decisions. The goal of this paper is to
understand the economic value of future market information, as energy
storages mitigate market uncertainty by forward-looking strategies. In
particular, we address the following research questions:

\begin{itemize}
\item From storages' perspective, what are the optimal decentralized buying
or selling quantities, when the energy prices are both uncertain and
variable over time?

\item At a system level, what is a good information management policy, to
coordinate storages' actions and improve their profitability?
\end{itemize}

To be specific, we provide a stylized model of optimal storage planning
under private market price forecasting. Decentralized decision makers have
to consider their optimal strategies in a competitive environment with
strategic interactions. In terms of information management, we consider the
following possible policy interventions: (1) providing a publicly available
market forecasting channel; (2) encouraging decentralized storages to share
their private forecasts with each other; (3) release additional market
information to a targeted subset of storages exclusively.

The rest of this paper is organized as follows. Section~\ref{s-lit} reviews
relevant literature. Section~\ref{s-model} introduces our model setup. In
Section \ref{s-analysis}, we carry out the analysis for two basic
(simplified) models. In Section \ref{s-policy}, we describe several policy
intervention solutions. In Section \ref{s-extend}, we extend the basic
models in several directions. Section~\ref{s-con} concludes this paper with
a discussion of the future research directions.

\section{Literature Review}

\label{s-lit} Our work contributes to the literature on oligopoly energy
market. The Cournot setup is a good approximation to some energy markets,
e.g., California's electricity industry as has been demonstrated in %
\citep{borenstein1999empirical}. Similar empirical work has been done in New
Zealand's electricity markets \citep{scott1996modelling}. This paper is
partly inspired by \cite{acemoglu2015competition}, wherein they also
consider a competitive energy market with highly asymmetric information
structures. While they seek to mitigate uncertainty and economic
inefficiency via contractual designs, we aim the same target with
informational interventions. We also focus on energy storages wherein they
consider energy producers.

In terms of energy storage modeling, our model extends a similar work
presented in \cite{contreras2015cooperation}, wherein they assume complete
information and deterministic demand function. The fundamental inefficiency
of such an energy market is driven by highly volatile local market
conditions (e.g. electricity prices), for instance, due to intermittency in
the renewable energy supply. For this reason, there is a growing literature
on the use of an energy storage system to improve integration of the
renewable energy \citep{Dicorato2012,Shu2014}. With this motivation (while
abstracting away from the physical characteristics of the renewables), our
model is closely related to this literature by incorporating both
intertemporal variability and uncertainty (exogenous market price shocks).
Given such an environment, storages will be foresighted and joint storage
planning and forecasting have been reported consequently %
\citep{Li2015,Haessig2015}. While this literature is mostly
simulation-based, our model admits tractable analysis and interpretable
structural results.

Beyond discussions on distributed storage planning and control, we put
emphasis on information management at a system level. In a deregulated
environment, it is natural to consider that the competing storages do not
observe each other's private information (private energy price forecast).
Therefore, the storages have to estimate each other's private forecast and
conjecture on how each other's action depends on its forecast. This
strategic interaction poses technical challenge, which is new to the energy
literature. A similar problem is studied in \cite{Kamalina2014} in a
different context (generation capacity expansion), wherein no structural
results are available. \cite{Shahidehpour2005} and \cite{Langary2014} touch
upon this topic in the context of generating companies' supply function
equilibrium but resort to simulation. Furthermore, the private forecast in
our model is sequentially revealed at every periods, while the private
information in aforementioned literature is \emph{static} (viewed as
generation companies' attribute or type).

Finally, there is a long stream of literature in economics and operations
research in terms of information management in such a decentralized setting.
We consider a class of equilibria wherein each agent's action depends
linearly on its forecast and forms Bayes' estimator for others' forecasts.
The uniqueness of this equilibrium prediction is guaranteed by \cite%
{Radner1962}. The value of a public forecast in coordinating agents' actions
is pioneered in \cite{Morris2002}. The incentives for information sharing
are studied in \cite{Gal-Or1985}. A recent work has demonstrated the power
of targeted information release \cite{Zhou2016}. We systematically examine
those ideas in the context of distributed energy storage market.

\section{Model}

\label{s-model} \textbf{Market structure.} Consider $n$ storages who
purchase and sell substitutable energy through a common market. The
storages, indexed by a set $I=\{1,2,\cdots ,n\}$, are homogeneous \textit{ex
ante} and engage in a Cournot competition. Let $d_{i}^{[t]}$ denote the
energy purchased (when $d_{i}^{[t]}<0$) or sold (when $d_{i}^{[t]}>0$) by
the $i^{th}$ storage at time $t$, and the aggregate storage quantity is
denoted by $D^{[t]}=\sum_{i=1}^{i=n}d_{i}^{[t]}$. We model the demand side
by assuming that the actual market clearing price $P^{[t]}(D^{[t]})$ is
linear in $D^{[t]}$, i.e.,

\begin{equation}
P^{[t]}(D^{[t]})=\beta ^{\lbrack t]}-\gamma ^{\lbrack t]}D^{[t]}+\eta
^{\lbrack t]},\forall t\in T,
\end{equation}%
where a random variable $\eta ^{\lbrack t]}$ captures the market
uncertainty. $\beta ^{\lbrack t]}>0$ corresponds to market potential, which
also captures market variability since $\beta ^{\lbrack t]}$ is changing
over time. $\gamma ^{\lbrack t]}>0$ (price elasticity) captures the fact
that the market price decreases when the aggregate energy sold $D^{[t]}$
increases, as the market supply of energy increases.

To model storages' strategic interactions, our demand side setup corresponds
to a scenario wherein storages are not price-takers but enjoy market power.
This is supported by empirical evidence that energy prices vary in response
to loads generation, especially when the storages are of sufficient scale %
\citep{Sioshansi2009}. Even if the storages are small-scaled, economics
literature suggests that infinitesimal agents also act \emph{as if} they are
expecting the price-supply relationship \citep{Osborne2005}. Similar
assumptions on storages' market power and price-anticipatory behavior are
not uncommon in the literature \citep{Sioshansi2010}. Finally, with the
integration of renewables and consequently the ``merit order effect'', the
supply of renewable energy can drastically impact the spot equilibrium price
considering its negligible marginal costs.

We assume that the market uncertainty follows an autoregressive process such
that
\begin{equation}
\eta ^{\lbrack t+1]}=\delta \eta ^{\lbrack t]}+\epsilon _{t},\forall t\in
T,\eta ^{\lbrack 1]}\sim N(0,\alpha ^{-1}).
\end{equation}%
The parameter $\alpha $ is the initial information precision concerning the
market uncertainty \textit{a priori}. Standard assumptions for
autoregressive process require that $\left\vert \delta \right\vert <1$,
where $\epsilon _{t}$ is exogenous shock $\epsilon _{t}\sim N(0,\zeta ^{-1})
$. We choose this stochastic process as it is among the simplest ones which
capture intertemporal correlation while remaining realistic.

\textbf{Storage model}. Storages are agents who can buy energy at a certain
time period and sell it at another. The net energy purchased and sold across
time is required to be zero for every storage, i.e., $\sum_{t\in
T}d_{i}^{[t]}=0$, $\forall i\in I$. As we seek to emphasize the interaction
between storages, we also abstract away from other operational constraints
such as energy and/or power limits; Instead, they are modeled by a cost
function associated with each storage. The battery degradation, efficiency,
and/or energy transaction costs of storage $i$ are represented by the cost
function $c_{i}(\cdot )$. This treatment is similar to that in \cite%
{contreras2015cooperation}.

We assume that $c_{i}(d)=\varepsilon ^{\lbrack t]}\cdot d^{2}$ in the basic
models. More realistically, $c_{i}(d)$ will be a power function $\varepsilon
^{\lbrack t]}\cdot d^{x}$ wherein $x\in (1,2)$, as it is known that as the
depth of discharge increases, the costs of utilizing storage increase faster
than linear. We can show that most of our results remain robust when the
cost function $c_{i}(d)$ is generalized within this region. We choose $x=2$
for a clear presentation of results. We also generalize the basic model to
consider heterogeneous cost functions in the extensions. To summarize this
discussion, the payoff of storage $i$ can be expressed as

\begin{equation}
\pi _{i}\left( d_{i}^{[t]},t\in T\right) =\sum_{t\in T}\left[
P^{[t]}(D^{[t]})\cdot d_{i}^{[t]}-\varepsilon ^{\lbrack t]}\cdot \left(
d_{i}^{[t]}\right) ^{2}\right] .
\end{equation}

\textbf{Information structure and sequence of events}. Storage $i$ has a
private forecast channel for the market condition. At the beginning of
period $t$, storage $i$ receives a private forecast $x_{i}^{[t]}$ with
precision $\rho $, i.e., $x_{i}^{[t]}=\eta ^{\lbrack t]}+\xi _{i}^{[t]}$,
where $\xi _{i}^{[t]}\sim $ $N(0,\rho ^{-1})$, for $\forall i\in I$. The
realizations of the forecasts are private, while their precision is common
knowledge. In this paper, we use ``forecast" and ``information"
interchangeably depending on the context.

At time $t$, the sequence of events proceeds as follows: (1) The storages
observe (the realizations of) their private forecasts; (2) Each storage
decides the purchase or selling quantities based on their information,
anticipating the rational decisions of the other storages; (3) The actual
market price is realized and the market is cleared for period $t$.
\begin{table}[h!]
\caption{Summary of nomenclature.}%
\begin{tabular}{ll}
\hline
$I$ & The set of energy storages. $|I|=n$. \\
$J$ & The set of targeted information release recipient. $|J|=m$. \\
$T$ & The set of time periods. $|T|=L$. \\
$X_{i}^{[t]}$ & Information set of storage $i$ in period $t$. \\ \hline
$A$ & Equilibrium base storage quantity. \\
$C$ & Equilibrium response factors with respect to private forecasts. \\
$B$ & Equilibrium response factor towards the public forecast. \\ \hline
$d_{i}^{[t]}$ & The amount of energy purchased or sold by the $i^{th}$
storage at time $t$. \\
$D^{[t]}$ & The aggregate storage quantity $D^{[t]}=%
\sum_{i=1}^{i=n}d_{i}^{[t]} $. \\
$P^{[t]}$ & Market clearing energy price at time $t$. \\
$x_{i}^{[t]}$ & Private forecast received by storage $i$ regarding market
uncertainty. \\
$x_{0}^{[t]}$ & Public forecast. \\ \hline
$\eta ^{\lbrack t]}$ & Market price uncertainty at time $t$. \\
$\beta ^{\lbrack t]}$ & Market potential. \\
$\gamma ^{\lbrack t]}$ & Energy price elasticity. \\
$\alpha $ & Information precision for market uncertainty \textit{a priori}.
\\
$\delta$ & Autoregression parameter for market uncertainty. \\
$\epsilon _{t}$ & Exogeneous price shock at time $t$. \\
$\zeta$ & Precision of $\epsilon _{t}$. \\
$\varepsilon ^{\lbrack t]}$ & Quadratic energy storage cost coefficient. \\
$\rho$ & Precision of the private forecast $x_{i}^{[t]}$. \\
$\pi _{i}$ & Storage $i$'s payoff function. \\
$\sigma$ & Precision of the public forecast $\xi_{0}^{[t]}$. \\
$\xi _{i}^{[t]}$ & Noise of the private information channel by storage $i$.
\\
$\xi_{0}^{[t]}$ & Noise of the public information channel. \\ \hline
\end{tabular}%
\end{table}

\section{Model Analysis}

\label{s-analysis}

\subsection{Centralized Energy Storage Model}

We begin by considering a single storage and thus dropping the subscript in
this section. The optimal storage quantities are obtained by solving
\begin{equation*}
\max_{d^{[t]},t\in T}\sum_{t\in T}\mathbb{E}\left[ \left.
\begin{array}{c}
P^{[t]}\cdot d^{[t]} \\
-\varepsilon ^{\lbrack t]}\cdot \left( d^{[t]}\right) ^{2}%
\end{array}%
\right\vert X^{[t]}\right] ,
\end{equation*}%
subject to
\begin{equation}
\mathbb{E}\left[ \left. \sum_{t\in T}d^{[t]}\right\vert X^{[t]}\right] =0,
\end{equation}%
for any sample path generated by $\{X^{[t]}\}^{\prime }s$, wherein $%
X^{[t]}=\left\{ x^{[1]},\eta ^{\lbrack 1]},...,\eta ^{\lbrack
t-1]},x^{[t]}\right\} $ indicates the corresponding information set. For
clarity of presentation, we solve for an optimal $d^{[t]}$ in an arbitrary
period $t$. In addition, we drop the superscript for $\gamma $ and $%
\varepsilon $ to focus on the intertemporal variability solely in market
price.

Notice that in period $t$, $d^{[t+1]},...,d^{[L]}$ will be anticipated
future optimal quantities based on the current information set $X^{[t]}$,
whereas $d^{[1]},...,d^{[t-1]}$ will be previous decisions realized to the
storage. To avoid confusion, we denote their solutions by a general $\mathbb{%
E}_{t}d^{[\tau ]}$, $\tau =1,...,L$. Under this notation, $\mathbb{E}%
_{t}d^{[\tau ]}$, $\tau =1,...,t-1$ will be known data, $\mathbb{E}%
_{t}d^{[t]}$ is the decision to be made in period $t$, and $\mathbb{E}%
_{t}d^{[\tau ]}$, $\tau =t+1,...,L$ will be the anticipated future optimal
quantities. It should be emphasized that $\mathbb{E}_{t}d^{[\tau ]}$ may not
be the same as the actual decision made in a future period $\tau $, for $%
\tau =t+1,...,L$. The timeline in this model is shown in Figure \ref%
{fig:timelinecentral}.

\begin{figure}[!h]
\centering
\advance\leftskip-4cm \includegraphics[scale=0.9]{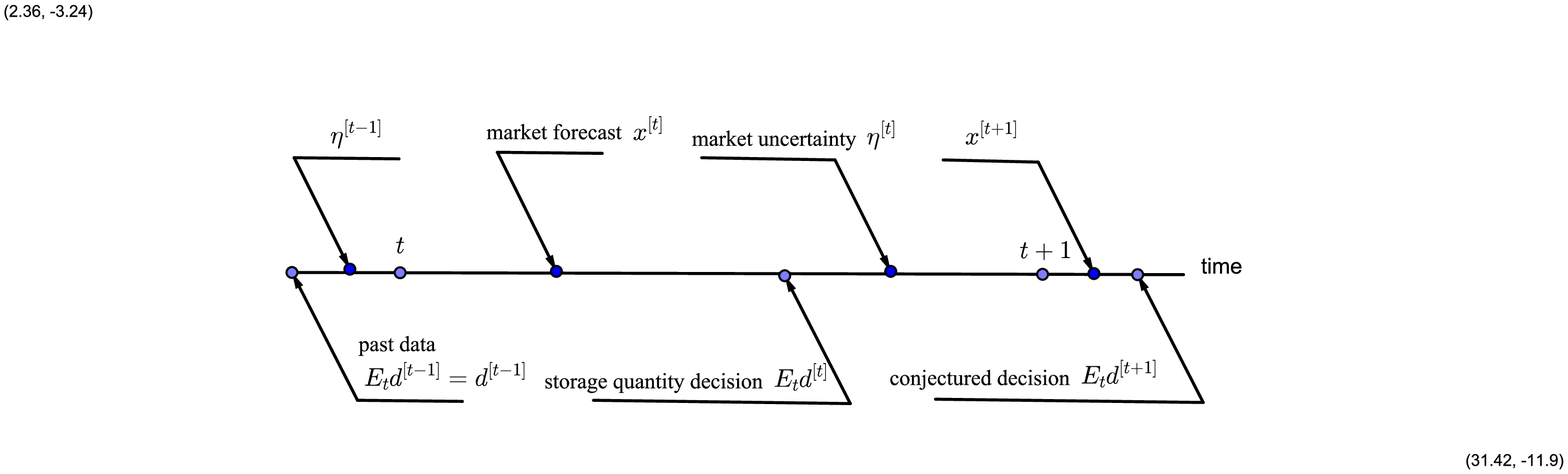}
\caption{Timeline of the centralized storage model.}
\label{fig:timelinecentral}
\end{figure}

By the \textit{Principal of Optimality}, in period $t$, we use the following
induced sub-problem to find $\mathbb{E}_{t}d^{[t]}$ (optimal solution of $%
d^{[t]}$):

\begin{equation*}
\max_{d^{[t]},d^{[t+1]},...,d^{[L]}}\sum_{\tau =t}^{\tau =L}\mathbb{E}\left[
\left. \left( \beta ^{\lbrack \tau ]}-\gamma d^{[\tau ]}+\eta ^{\lbrack \tau
]}\right) d^{[\tau ]}-\varepsilon \cdot \left( d^{[\tau ]}\right)
^{2}\right\vert X^{[t]}\right] ,
\end{equation*}%
subject to
\begin{equation}
\mathbb{E}\left[ \left. \sum_{\tau =t}^{\tau =L}d^{[\tau ]}\right\vert
X^{[t]}\right] =-\sum_{\tau =1}^{\tau =t-1}d^{[\tau ]}.
\end{equation}

\begin{pro}
\label{pro_central} The optimal storage quantity in period $t$ is denoted by
\begin{eqnarray}
\mathbb{E}_{t}d^{[t]} &=&\underset{\text{base storage quantity}}{\underbrace{%
\frac{\beta ^{\lbrack t]}-\frac{\sum_{\tau =t}^{\tau =L}\beta ^{\lbrack \tau
]}}{L-t+1}}{2\left( \varepsilon +\gamma \right) }-\frac{\sum_{\tau =1}^{\tau
=t-1}d^{[\tau ]}}{L-t+1}}} \\
&&+\underset{\text{response factor}}{\underbrace{\frac{\left( 1-\frac{%
\sum_{\tau =t}^{\tau =L}\delta ^{\tau -t}}{L-t+1}\right) }{2\left(
\varepsilon +\gamma \right) }}}\cdot \mathbb{E}\left[ \left. \eta ^{\lbrack
t]}\right\vert X^{[t]}\right],
\end{eqnarray}%
for $t=1,2,...,L-1$, wherein $\mathbb{E}\left[ \left. \eta ^{\lbrack
t]}\right\vert X^{[t]}\right] =\frac{\rho }{\rho +\zeta }x^{[t]}+\frac{\zeta
\delta }{\rho +\zeta }\eta ^{\lbrack t-1]}$, for $t=2,...,L$, and $\mathbb{E}%
\left[ \left. \eta ^{\lbrack 1]}\right\vert X^{[1]}\right] =\frac{\rho }{%
\rho +\alpha }x^{[1]}$. In the final period, $\mathbb{E}_{L}d^{[L]}=-\sum_{%
\tau =1}^{\tau =L-1}d^{[\tau ]}$.
\end{pro}

In this proposition, we are able to derive the optimal storage quantities in
closed form, comprising of two parts: The first part is the \textit{base
storage quantity}, and the second part is the \textit{response factor}
multiplied by an estimation of the market price uncertainty. The base
storage quantity decreases in $\varepsilon $ (cost coefficients with respect
to the depth of discharge) and $\gamma $ ( energy price elasticity). The
market price uncertainty is estimated by a convex combination of a current
forecast and a last-period observation: $\mathbb{E}\left[ \left. \eta
^{\lbrack t]}\right\vert X^{[t]}\right] =\frac{\rho }{\rho +\zeta }x^{[t]}+%
\frac{\zeta \delta }{\rho +\zeta }\eta ^{\lbrack t-1]}$. The weighting
factors are proportional to their relative precision levels $\rho $ and $%
\zeta $. In additional, the last-period observation weights more when the
intertemporal correlation is stronger ($\delta $ is higher).

From this proposition, we can clearly see that the optimal storage
quantities depend on both variability (captured by market potential
coefficient $\beta ^{\lbrack t]}$) and uncertainty (captured by $\eta
^{\lbrack t]}$) in market price. The base storage quantity demonstrates a
downward distortion to the per-stage optimal storage quantity $\frac{\beta
^{\lbrack t]}}{2\left( \varepsilon +\gamma \right) }$: The component $\frac{%
\sum_{\tau =t}^{\tau =L}\frac{\beta ^{\lbrack \tau ]}}{L-t+1}}{2\left(
\varepsilon +\gamma \right) }$ is the average per-stage optimal storage
quantity across all future periods, and the component $\frac{\sum_{\tau
=1}^{\tau =t-1}d^{[\tau ]}}{L-t+1}$ is subtracted to compensate the existing
energy storage level built up in the past. The downward distortion ensures
that the overall storage quantities offset each other, i.e., $\mathbb{E}%
\left[ \left. \sum_{t\in T}d^{[t]}\right\vert X^{[t]}\right] =0$.

\subsection{Decentralized Two-Period Model}

In this section, we recover the superscript for $\gamma $ and $\varepsilon $%
. We also recover the subscript for $x_{i}^{[t]}$ since the storages observe
heterogeneous information. To characterize the equilibrium outcome under
such highly asymmetric information structure, we first introduce our
solution concept.

\textbf{Equilibrium concept}. Storage $i$ chooses a storage quantity $%
d_{i}^{[1]}$ to maximize $\mathbb{E}[\pi _{i}|x_{i}^{[1]}]$, by forming an
expectation of the other producers' production levels $\mathbb{E}%
(d_{j}^{[1]}|x_{i}^{[1]})$, for $\forall j\neq i$. $d_{i}^{[2]}$ is
determined thereafter, due to the constraint $d_{i}^{[1]}+d_{i}^{[2]}=0$, $%
\forall i\in I$. We focus exclusively on the \textit{linear symmetric
Bayesian-Nash equilibrium}, i.e., $d_{i}^{[1]}=A+Cx_{i}^{[1]},$ for some
constants $A$ and $C$. We can interpret $A$ as the \textit{base storage
quantity}; $C$ as the \textit{response factors} with respect to the forecast
$x_{i}^{[1]}$, respectively.

\begin{pro}
\label{prop_two period} For a two-period model under private market
forecasting, the storage quantity in the linear symmetric Bayesian-Nash
equilibrium is $d_{i}^{[1]}=A+Cx_{i}^{[1]}$ for every storage, wherein
\begin{equation*}
A=\frac{\beta ^{\lbrack 1]}-\beta ^{\lbrack 2]}}{2\left( \varepsilon
^{\lbrack 1]}+\varepsilon ^{\lbrack 2]}\right) +(n+1)\left( \gamma ^{\lbrack
1]}+\gamma ^{\lbrack 2]}\right) },
\end{equation*}

\begin{equation*}
C=\frac{(1-\delta )\rho }{\left[
\begin{array}{c}
(n-1)\left( \gamma ^{\lbrack 1]}+\gamma ^{\lbrack 2]}\right) \rho + \\
2\left( \varepsilon ^{\lbrack 1]}+\varepsilon ^{\lbrack 2]}+\gamma ^{\lbrack
1]}+\gamma ^{\lbrack 2]}\right) (\alpha +\rho )%
\end{array}%
\right] }.
\end{equation*}
\end{pro}

In period $t=1$, the base storage quantity $A$ is positive (selling energy)
if and only if $\beta ^{\lbrack 1]}>\beta ^{\lbrack 2]}$, i.e., the energy
price decreases at $t=2$. Conversely, the storage buys energy ($A<0$) when
it can be sold at a higher price at $t=2$ ($\beta ^{\lbrack 1]}<\beta
^{\lbrack 2]}$). Consistent with our results in the single storage model,
the base selling quantity decreases in $\sum_{t\in T}\varepsilon ^{\lbrack
t]}$ (cost coefficients with respect to the depth of discharge), and in $%
\sum_{t\in T}\gamma ^{\lbrack t]}$ (aggregate energy price elasticity).
Furthermore, when $A>0$, it decreases in the number of storages. This is
because that the market competition is more intense as the number of Cournot
competitors increases, and consequently, the price of energy decreases. The
converse is true when $A<0$.

The reaction to private forecast $C$ is more aggressive when its precision $%
\rho $ increases, as a storage relies more on an accurate market forecast.
The reactions to forecasts are less aggressive when the intertemporal
correlation $\delta $ increases. In the extreme case where $\delta =1$, the
storage does not respond to forecasts. Intuitively, this is because that any
action (either buy or sell) in response to a market forecasts at $t=1$ will
be offset by a reverse operation (under the energy balance constraints $%
\sum_{t\in T}d_{i}^{[t]}=0$) at $t=2$ when the market condition remains the
same.

\begin{pro}
\label{prop:twoperiodpayoff} When there are a large number of storages,
every storage's equilibrium payoff converges asymptotically to
\begin{equation}
\lim_{n\rightarrow \infty }\mathbb{E}[\pi _{i}]\rightarrow
\sum_{t=1,2}\left( \varepsilon ^{\lbrack t]}+\gamma ^{\lbrack t]}\right) %
\left[
\begin{array}{c}
\left( \frac{\beta ^{\lbrack 1]}-\beta ^{\lbrack 2]}}{\gamma ^{\lbrack
1]}+\gamma ^{\lbrack 2]}}\right) ^{2}+ \\
\left( \frac{1-\delta }{\gamma ^{\lbrack 1]}+\gamma ^{\lbrack 2]}}\right)
^{2}\rho ^{-1}%
\end{array}%
\right] \cdot n^{-2}.
\end{equation}
\end{pro}

When there are a large number of storages, $\lim_{n\rightarrow \infty }%
\mathbb{E}[\pi _{i}]$ decreases in $\rho $. This result demonstrates the
negative economic value of a private forecast, and is interpreted as \textit{%
competition effect}. Although storages' private forecasts are independent,
they lead to similar reaction to a market price shock. When the number of
storages is large, the purchase or selling quantity responses are
exaggerated and the over-precision in forecasts can lead to even lower
payoffs. The inverse-square decay with respect to the number of storages
demonstrates the impact of competition intensity.

\begin{figure}[]
\begin{center}
\subfigure[Impact of autoregression parameter $\delta$.]{            \label{fig:first}
            \includegraphics[width=0.5\textwidth]{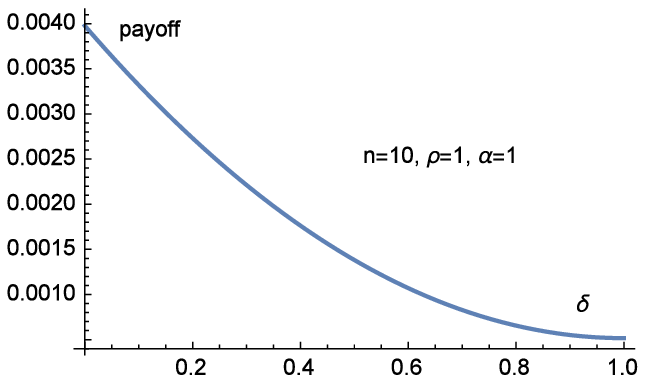}
        }
\subfigure[Impact of market uncertainty parameter $\alpha$.]{           \label{fig:second}
           \includegraphics[width=0.5\textwidth]{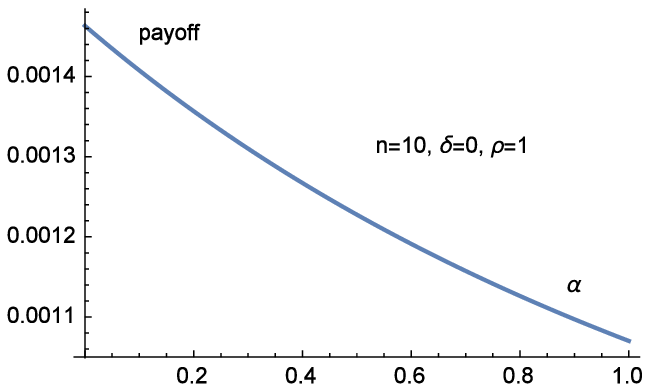}
        }\\[0pt]
\subfigure[Impact of the number of storages $n$.]{            \label{fig:third}
            \includegraphics[width=0.5\textwidth]{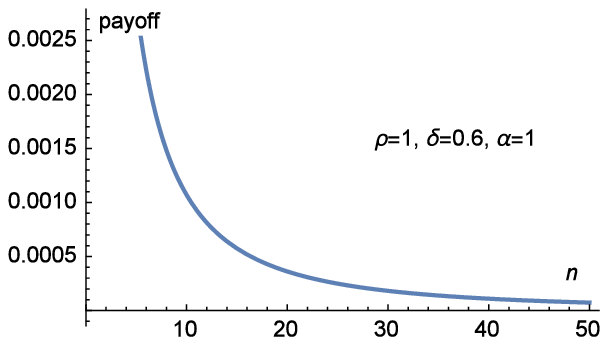}
        }
\subfigure[Impact of private information precision $\rho$.]{            \label{fig:fourth}
            \includegraphics[width=0.5\textwidth]{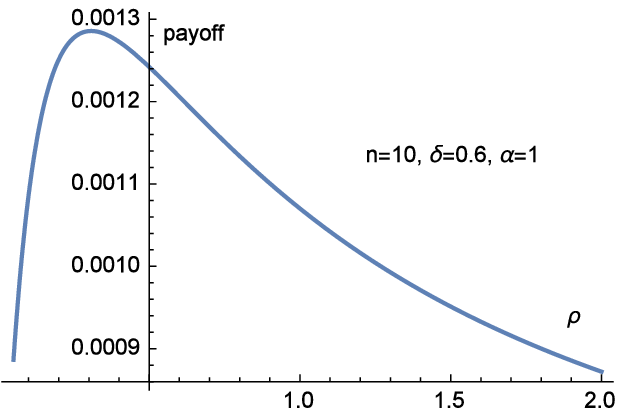}
        }
\end{center}
\caption{ Sensitivity analysis of individual storage's payoff with respect
to model primitives. }
\label{fig:sensitivity4}
\end{figure}

We complement the analysis with the following numerical example, wherein $%
\beta ^{\lbrack 1]}-\beta ^{\lbrack 2]}=1$, $\varepsilon ^{\lbrack
1]}=\varepsilon ^{\lbrack 2]}=\gamma ^{\lbrack 1]}=\gamma ^{\lbrack 2]}=1$,
and the payoff is calculated with a finite number of storages. We summarize
the results in Figure \ref{fig:sensitivity4}. A storage's payoff decreases
in the autoregressive parameter $\delta$. As a storage's reactions to
forecasts are less responsive when the intertemporal correlation ($\delta$)
increases, the economic value of market information also decreases. For a
similar reason, a storage's payoff also decreases in market uncertainty
parameter $\alpha$, as the storage's reactions to forecasts are also less
responsive in this case. In addition, the inverse-square decay with respect to
the number of storages confirms our result in the asymptotic analysis.
Finally, a storage's payoff first increases then decreases in its private
information precision. When private information is scarce, the economic
value of private information increases in its precision as it mitigate
uncertainty. However, when the information precision further increases, the
competition effect dominates and the payoff decreases.

\section{Operational Policy Analysis}

\label{s-policy}

\subsection{Public Forecast Provision}

Suppose that instead of private forecasts, all the storages receive a public
forecast $x_{0}^{[t]}=\eta ^{\lbrack t]}+\xi _{0}^{[t]}$, where $\xi
_{0}^{[t]}\sim $ $N(0,\sigma ^{-1})$. In this section, we analyze the
possibility for a public forecast to coordinate storages' actions. This can
be potentially provided by the aggregator. Following a similar analysis as
in the private forecasting model, we assume that $d_{i}^{[1]}=A+Bx_{0}^{[1]}$%
, where $B$ is the response factor towards the public forecast.

\begin{pro}
\label{prop_public} For a two-period model under public market forecasting,
each storage's equilibrium storage quantity in the linear symmetric
Bayesian-Nash equilibrium is $d_{i}^{[1]}=A+Bx_{0}^{[1]}$, wherein
\begin{equation*}
A=\frac{\beta ^{\lbrack 1]}-\beta ^{\lbrack 2]}}{2\left( \varepsilon
^{\lbrack 1]}+\varepsilon ^{\lbrack 2]}\right) +(n+1)\left( \gamma ^{\lbrack
1]}+\gamma ^{\lbrack 2]}\right) },
\end{equation*}%
\begin{equation*}
B=\frac{(1-\delta )\sigma /(\alpha +\sigma )}{2\left( \varepsilon ^{\lbrack
1]}+\varepsilon ^{\lbrack 2]}\right) +(n+1)\left( \gamma ^{\lbrack
1]}+\gamma ^{\lbrack 2]}\right) }.
\end{equation*}%
The equilibrium payoff
\begin{equation*}
\lim_{n\rightarrow \infty }\mathbb{E}[\pi _{i}]\rightarrow
\sum_{t=1,2}\left( \varepsilon ^{\lbrack t]}+\gamma ^{\lbrack t]}\right) %
\left[
\begin{array}{c}
\left( \frac{\beta ^{\lbrack 1]}-\beta ^{\lbrack 2]}}{\gamma ^{\lbrack
1]}+\gamma ^{\lbrack 2]}}\right) ^{2} \\
+\left( \frac{1-\delta }{\gamma ^{\lbrack 1]}+\gamma ^{\lbrack 2]}}\right)
^{2}\cdot \frac{\sigma }{\left( \alpha +\sigma \right) ^{2}}%
\end{array}%
\right] \cdot n^{-2}.
\end{equation*}
\end{pro}

Similar to the private forecasting model, the reaction to public forecast $B$
is more aggressive when its precision $\sigma $ increases. $%
\lim_{n\rightarrow \infty }\mathbb{E}[\pi _{i}]$ is pseudo-concave in the
public forecast precision $\sigma $, and thus reaches global maximum when $%
\sigma =\alpha $. Notice that, when $\sigma >\alpha $, $\lim_{n\rightarrow
\infty }\mathbb{E}[\pi _{i}]$ decreases in $\sigma $, which demonstrate the
negative economic value of a public forecast. We interpret this as \textit{%
congestion effect}: Intuitively, when $\sigma >\alpha $, over-reaction to a
public forecast leads to either too much purchase quantity (when forecast is
favorable) or too much selling quantity (unfavorable forecast) from all
storages. The public forecast becomes a herding signal. It is beneficial to
maintain certain exclusiveness of a public forecast. For given information
provision, a storage's payoff suffers inverse-square decay in the number of
storages, as the economic value of a public forecast is diluted when more
storages respond to it.

\begin{figure}[]
\begin{center}
\subfigure[Impact of number of storages $n$.]{            \label{fig:pubfirst}
            \includegraphics[width=0.5\textwidth]{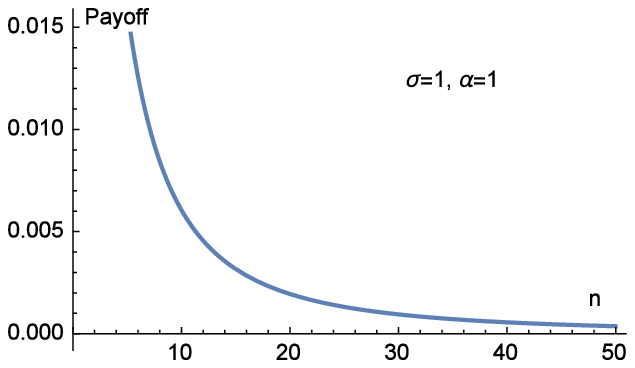}
        }
\subfigure[Impact of public information precision $\sigma$.]{           \label{fig:pubsecond}
           \includegraphics[width=0.5\textwidth]{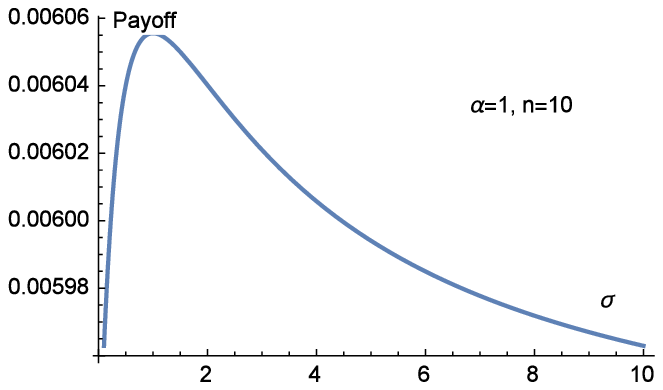}
        }
\end{center}
\caption{Sensitivity analysis of individual storage's payoff under public
information provision.}
\label{fig:sensitivity2}
\end{figure}

We complement the analysis with the following numerical example summarized
in Figure \ref{fig:sensitivity2}, wherein $\beta ^{\lbrack 1]}-\beta
^{\lbrack 2]}=1$, $\varepsilon ^{\lbrack 1]}=\varepsilon ^{\lbrack
2]}=\gamma ^{\lbrack 1]}=\gamma ^{\lbrack 2]}=1$, and the payoff is
calculated with a finite number of storages. Again, a storage's payoff
decreases to the inverse-square of the number of storages as shown in the
analytical result. A storage's payoff first increases (due to positive
economic value) then decreases (due to the congestion effect) in the
precision of the public market forecast.

\subsection{Encourage Information Sharing}

Suppose that the storages pool their private forecasts $x_{i}^{[1]}$
together. In this case, it can be checked that it is equivalent for them to
observe a public forecast $x_{0}^{[1]}$\ with precision $\sigma =n\rho $:

\begin{equation*}
\mathbb{E}[\eta ^{\lbrack 1]}|x_{1}^{[1]},...,x_{n}^{[1]}]=\frac{\rho }{%
\alpha +n\rho }\sum_{i\in I}x_{i}^{[1]},\mathbb{E}[\eta ^{\lbrack
1]}|x_{0}^{[1]}]=\frac{\sigma }{\alpha +\sigma }x_{0}^{[1]},
\end{equation*}%
and it can be checked that these two estimators are stochastically
equivalent, i.e., both $N(0,\frac{n\rho }{(\alpha +n\rho )^{2}})$.

\begin{figure}[]
\begin{center}
\subfigure[Small number of energy storages.]{            \label{fig:comfirst}
            \includegraphics[width=0.5\textwidth]{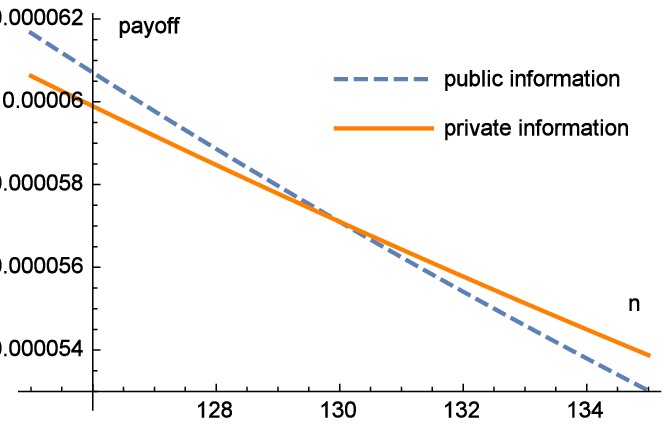}
        }
\subfigure[Large number of energy storages.]{           \label{fig:comsecond}
           \includegraphics[width=0.5\textwidth]{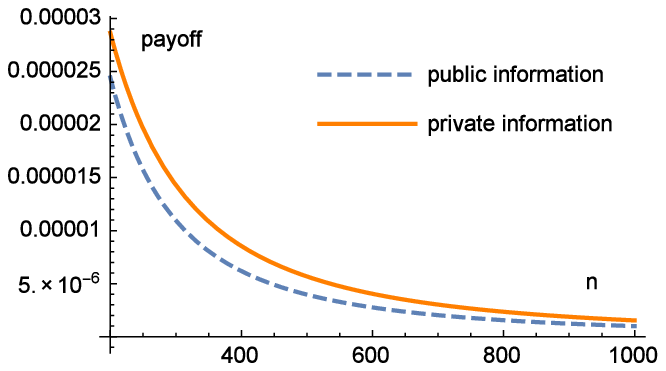}
        }
\end{center}
\caption{Individual payoff comparison with and without information sharing.}
\label{fig:sensitivity3}
\end{figure}

Therefore, we can calculate the corresponding payoffs under pooling private
forecasts:

\begin{eqnarray*}
\lim_{n\rightarrow \infty }\mathbb{E}[\pi _{i};\sigma &=&n\rho ]\rightarrow
\sum_{t=1,2}\left( \varepsilon ^{\lbrack t]}+\gamma ^{\lbrack t]}\right)
\cdot \left[
\begin{array}{c}
\left( \frac{\beta ^{\lbrack 1]}-\beta ^{\lbrack 2]}}{\gamma ^{\lbrack
1]}+\gamma ^{\lbrack 2]}}\right) ^{2} \\
+\left( \frac{1-\delta }{\gamma ^{\lbrack 1]}+\gamma ^{\lbrack 2]}}\right)
^{2}\cdot \left( n\rho \right) ^{-1}%
\end{array}%
\right] \cdot n^{-2} \\
&<&\sum_{t=1,2}\left( \varepsilon ^{\lbrack t]}+\gamma ^{\lbrack t]}\right)
\left[ \left( \frac{\beta ^{\lbrack 1]}-\beta ^{\lbrack 2]}}{\gamma
^{\lbrack 1]}+\gamma ^{\lbrack 2]}}\right) ^{2}+\left( \frac{1-\delta }{%
\gamma ^{\lbrack 1]}+\gamma ^{\lbrack 2]}}\right) ^{2}\rho ^{-1}\right]
\cdot n^{-2}.
\end{eqnarray*}%
By comparing this payoff with that under private forecasts, we find that the
economic value of forecast under information sharing is (to an order of
magnitude in the number of storages) lower than that under private
forecasts. There will be no incentive for the storage to share information
with each others. From this analysis, we find that communication among the
storages fails to achieve a coordinated effort to increase market
efficiency. To maintain the exclusiveness of their private forecasts, the
decentralized storages should not be encouraged to share market information
in this regime. This result is confirmed by a numerical analysis summarized
in Figure \ref{fig:sensitivity3} when there are a large number of storages.
However, when the number of storages is small, it is possible that every
storage is better off by information sharing. This regime is possible when
each private forecast is extremely fuzzy, and pooling them together can
amplify the market signal.

\subsection{Targeted Information Release}

Now that we know the exclusiveness of a market forecast is important, we
analyze an alternative policy intervention through public information
channel. Suppose that the aggregator/government offers a public forecast
only to a subset of storages $J$ ($|J|=m\leq n$). For informed storages
(ones who receive the public forecast), their $d_{i}^{[1]}=A+Bx_{0}^{[1]}$, $%
\forall i$ $\in J$. For uninformed storages, their $d_{i}^{[1]}=C$, $\forall
i$ $\in I-J$. $A$, $B$ and $C$ are all unknown constant coefficients.

\begin{pro}
\label{prop_target} For a two-period model under targeted information
release, storages' equilibrium storage quantities in the linear
Bayesian-Nash equilibrium are $d_{i}^{[1]}=A+Bx_{0}^{[1]}$, $\forall i$ $\in
J$, and $d_{i}^{[1]}=C$, $\forall i$ $\in I-J$, wherein
\begin{equation*}
A=C=\frac{\beta ^{\lbrack 1]}-\beta ^{\lbrack 2]}}{2\left( \varepsilon
^{\lbrack 1]}+\varepsilon ^{\lbrack 2]}\right) +(n+1)\left( \gamma ^{\lbrack
1]}+\gamma ^{\lbrack 2]}\right) },
\end{equation*}%
\begin{equation*}
B=\frac{(1-\delta )\sigma /(\alpha +\sigma )}{2\left( \varepsilon ^{\lbrack
1]}+\varepsilon ^{\lbrack 2]}\right) +(m+1)\left( \gamma ^{\lbrack
1]}+\gamma ^{\lbrack 2]}\right) }.
\end{equation*}%
The storages' aggregate payoff $\sum_{i\in I}\mathbb{E}[\pi _{i}]$ is
maximized when the population of information recipient%
\begin{equation}
m=1+\frac{2\left( \varepsilon ^{\lbrack 1]}+\varepsilon ^{\lbrack 2]}\right)
}{\gamma ^{\lbrack 1]}+\gamma ^{\lbrack 2]}}.
\end{equation}
\end{pro}

In this case, the storages' payoffs are stratified, due to their
asymmetrical informational status. The fact that an interior solution (of $m$%
) exists suggests a trade-off between the economic value of a public
forecast in coordinating the storages' actions, and the congestion effect
due to the lack of exclusiveness of such information dissemination.

\section{Model Generalizations}

\label{s-extend}

\subsection{Multi-Period Model}

In this section, we demonstrate that the model can be extended in multiple
directions. A multi-period version of this problem has to be solved
recursively using backward induction while unfolding the information set
throughout the process. Instead, we analyze a relaxed problem. In this case,
equilibrium characterization requires solving the following optimization
problems:

\begin{equation*}
\max_{d_{i}^{[t]},t\in T}\sum_{t\in T}\mathbb{E}\left[ \left.
\begin{array}{c}
P^{[t]}(D^{[t]})\cdot d_{i}^{[t]} \\
-\varepsilon ^{\lbrack t]}\cdot \left( d_{i}^{[t]}\right) ^{2}%
\end{array}%
\right\vert X_{i}^{[t]}\right] ,
\end{equation*}%
subject to
\begin{equation}
\mathbb{E}\left[ \sum_{t\in T}d_{i}^{[t]}\right] =0.
\end{equation}%
$X_{i}^{[t]}=\left\{ x_{0}^{[1]},x_{i}^{[1]},\eta ^{\lbrack 1]},...,\eta
^{\lbrack t-1]},x_{0}^{[t]},x_{i}^{[t]}\right\} $ indicates the
corresponding information set. Notice that we simultaneous incorporate
private forecasts $\{x_{i}^{[t]}\}$'s and a public forecast $x_{0}^{[t]}$.
The storage quantity will be $%
d_{i}^{[t]}=A^{[t]}+B^{[t]}x_{0}^{[t]}+C^{[t]}x_{i}^{[t]}$ for some unknown
coefficients $A^{[t]}$, $B^{[t]}$ and $C^{[t]}$. This is a relaxation
because an exact solution requires that $\mathbb{E}\left[ \left. \sum_{t\in
T}\left( d_{i}^{[t]}\right) ^{\ast }\right\vert X_{i}^{[t]}\right] =0$, for
any sample path generated by $\{X_{i}^{[t]}\}^{\prime }s$.

\begin{pro}
\label{prop_multiple period} For a multi-period model under both private and
public market forecasting, the storage quantities in the linear symmetric
Bayesian-Nash equilibrium can be approximated by $%
d_{i}^{[t]}=A^{[t]}+B^{[t]}x_{0}^{[t]}+C^{[t]}x_{i}^{[t]}$, wherein
\begin{equation*}
A^{[t]}=\frac{\beta ^{\lbrack t]}-\lambda }{2\varepsilon ^{\lbrack
t]}+(n+1)\gamma ^{\lbrack t]}},
\end{equation*}

\begin{equation*}
C^{[t]}=\frac{\rho }{\left[ 2\left( \varepsilon ^{\lbrack t]}+\gamma
^{\lbrack t]}\right) (\alpha +\sigma +\rho )+(n-1)\gamma ^{\lbrack t]}\rho %
\right] },
\end{equation*}

\begin{equation*}
B^{[t]}=\frac{\sigma -(n-1)\gamma ^{\lbrack t]}\sigma C}{\left[ (n+1)\gamma
^{\lbrack t]}+2\varepsilon ^{\lbrack t]}\right] (\alpha +\sigma +\rho )},
\end{equation*}%
and the Lagrangian multiplier
\begin{equation}
\lambda =\frac{\sum_{t\in T}\beta ^{\lbrack t]}\prod\limits_{\tau \neq t}%
\left[ 2\varepsilon ^{\lbrack \tau ]}+(n+1)\gamma ^{\lbrack \tau ]}\right] }{%
\sum_{t\in T}\prod\limits_{\tau \neq t}\left[ 2\varepsilon ^{\lbrack \tau
]}+(n+1)\gamma ^{\lbrack \tau ]}\right] }.
\end{equation}%
.
\end{pro}

The downside of this analysis is that we can not guarantee that $\mathbb{E}%
\left[ \left. \sum_{t\in T}\left( d_{i}^{[t]}\right) ^{\ast }\right\vert
X_{i}^{[t]}\right] =0$. Essentially, the storages reduce the baseline
quantity $A^{[t]}$ by their time-average so that the aggregate buy/sell
quantities sum up to zero in the statistical sense. This approximation of
the storages' actions ignores the intertemporal correlation of
uncertainties. Future research is needed for an exact analysis of the
full-fledged model.

\subsection{Heterogeneous Storages}

We model storages with heterogeneous physical attributes and information
status by assuming that the costs of utilizing storage $c_{i}(d)=\varepsilon
_{i}^{[t]}\cdot d^{2}$, and storage $i$ receives a private forecast $%
x_{i}^{[t]}$ with precision $\rho _{i}$. To illustrate the major points, we
extend the two-period model.

\begin{pro}
\label{prop_Hetero} For a two-period model under both private and public
market forecasting, the heterogeneous storage quantities in the linear
Bayesian-Nash equilibrium take the form of $%
d_{i}^{[1]}=A_{i}+B_{i}x_{0}^{[1]}+C_{i}x_{i}^{[1]}$, wherein $A_{i}$, $%
B_{i} $, and $C_{i}$ are given in the Appendix.
\end{pro}

As in the basic model with homogeneous storages, each storage holds their
own forecast but with varying precision:
\begin{equation}
\mathbb{E}[\eta ^{\lbrack 1]}|x_{0}^{[1]},x_{i}^{[1]}]=\frac{\sigma }{\alpha
+\sigma +\rho }x_{0}^{[1]}+\frac{\rho _{i}}{\alpha +\sigma +\rho _{i}}%
x_{i}^{[1]}.
\end{equation}%
The interesting new feature in this extension is that each storage ($i$)
needs to guess another storage's ($j$) quantity decision via its own
information set:

\begin{equation}
\mathbb{E}[d_{j}^{[1]}|x_{0}^{[1]},x_{i}^{[1]}]=A_{j}+B_{j}x_{0}^{[1]}+C_{j}%
\mathbb{E}[x_{j}^{[1]}|x_{0}^{[1]},x_{i}^{[1]}].
\end{equation}%
To obtain a correct conjecture, storage $i$ needs to estimate storage $j$'s
private forecast:

\begin{equation}
\mathbb{E}[x_{j}^{[1]}|x_{0}^{[1]},x_{i}^{[1]}]=\mathbb{E}[\eta ^{\lbrack
1]}+\xi _{j}^{[1]}|x_{0}^{[1]},x_{i}^{[1]}]=\mathbb{E}[\eta ^{\lbrack
1]}|x_{0}^{[1]},x_{i}^{[1]}].
\end{equation}%
Following a similar procedure as in the basic homogeneous model, we can
obtain the coefficients summarized in the Appendix. The dependency of the
value of both public and private forecasts on their precisions $\sigma $ and
$\rho _{i}$ is highly nonlinear. Due to the complicated payoff functional
forms of this general model, we start with the homogeneous model for a clear
presentation of results. It can be checked that some of our findings and
intuitions remain robust under this generalization. For example, the base
storage quantity $A$ is positive (selling energy) if and only if $\beta
^{\lbrack 1]}>\beta ^{\lbrack 2]}$, i.e., the energy price decreases at $t=2$%
. The value of a\ private forecast ($C_{i}^{2}\mathbb{E}[x_{i}^{[1]}]^{2}$)
is proportional to $\left[ \frac{(1-\delta )}{\alpha +\sigma +\rho _{i}}%
\right] ^{2}\rho _{i}$, and thus decreasing in the intertemporal correlation
$\delta $, as the reactions to forecasts are less aggressive when $\delta $
increases.

\section{Conclusion}

\label{s-con} In this paper, we propose stylized models of decentralized
energy storage planning under private and public market forecasting, when
energy prices are both uncertain and variable over time. We derive the
optimal buying or selling quantities for storages in a competitive
environment with strategic interactions. Coarsely speaking, a foresighted
storage will plan to buy energy when its price is low and sell when the
price high. The value of a private forecast decreases in the intertemporal
correlation of market price shock. We demonstrate the potentially negative
economic value of a private forecast, due to \emph{competition effect}: When
there are a large number of storages, the purchase or selling quantity
responses are exaggerated and the over-precision in forecasts can lead to
even lower payoffs. These fundamental observations are robust when we
generalize the model to multi-period or heterogeneous storages.

We also examine several information management policies to coordinate
storages' actions and improve their profitability. Firstly, we demonstrate
the potential negative economic value of a public forecast, due to \emph{%
congestion effect}: A precise public forecast leads to herding behavior, and
over-reaction to a public forecast leads to either too much purchase
quantity (when forecast is favorable) or too much selling quantity
(unfavorable forecast) from all storages. Secondly, we find that
communication among the storages could fail to achieve a coordinated effort
to increase market efficiency. The decentralized storages will not
participate in any information sharing program when there are a large number
of storages. Thirdly, we find it optimal to release additional information
to a subset of energy storages exclusively by targeted information release.

Future research is needed for a full-fledged analysis of a multi-period,
decentralized, and heterogeneous model. Another direction to go is to
incorporate operational constraints such as energy and/or power limits.
Explicit modeling of renewable energy generation will contribute to a
holistic understanding of the entire integrated system. Finally, information
management research in other energy markets is likely to be promising.

% or


\bibliographystyle{chicagoa}
\bibliography{InfoProvision}

\appendix

\section{Appendix. Proofs.}

\textbf{Proof of Proposition \ref{pro_central}}. Introduce a Lagrangian
multiplier $\lambda ^{\lbrack t]}$ to relax the conservation constraint. We
denote the Lagrangian by

\begin{eqnarray}
L &=&\sum_{\tau =t}^{\tau =L}l^{[\tau ]}  \notag \\
&=&\sum_{\tau =t}^{\tau =L}\mathbb{E}\left[ \left.
\begin{array}{c}
\left( \beta ^{\lbrack \tau ]}-\gamma d^{[\tau ]}+\eta ^{\lbrack \tau
]}\right) d^{[\tau ]} \\
-\varepsilon \cdot \left( d^{[\tau ]}\right) ^{2}-\lambda ^{\lbrack
t]}d^{[\tau ]}%
\end{array}%
\right\vert X^{[t]}\right] -\lambda ^{\lbrack t]}\sum_{\tau =1}^{\tau
=t-1}d^{[\tau ]}
\end{eqnarray}%
By the first-order condition $\frac{\partial l^{[\tau ]}}{\partial d^{[\tau
]}}=0$, for $\tau =t,...,L$ $\Rightarrow $

\begin{equation}
\mathbb{E}_{t}d^{[\tau ]}=\frac{\beta ^{\lbrack \tau ]}-\lambda ^{\lbrack
t]}+\mathbb{E}\left[ \left. \eta ^{\lbrack \tau ]}\right\vert X^{[t]}\right]
}{2\left( \varepsilon +\gamma \right) }.
\end{equation}%
Notice that $\mathbb{E}\left[ \left. \eta ^{\lbrack \tau ]}\right\vert
X^{[t]}\right] =\delta ^{\tau -t}\mathbb{E}\left[ \left. \eta ^{\lbrack
t]}\right\vert X^{[t]}\right] $, for $\tau =t,...,L$, and $\mathbb{E}\left[
\left. \sum_{\tau =t}^{\tau =L}d^{[\tau ]}\right\vert X^{[t]}\right]
=-\sum_{\tau =1}^{\tau =t-1}d^{[\tau ]}\Rightarrow $

\begin{equation}
\lambda ^{\lbrack t]}=\frac{\sum_{\tau =t}^{\tau =L}\beta ^{\lbrack \tau
]}+\sum_{\tau =t}^{\tau =L}\delta ^{\tau -t}\mathbb{E}\left[ \left. \eta
^{\lbrack t]}\right\vert X^{[t]}\right] +2\left( \varepsilon +\gamma \right)
\sum_{\tau =1}^{\tau =t-1}d^{[\tau ]}}{L-t+1}
\end{equation}%
Therefore,%
\begin{equation}
\mathbb{E}_{t}d^{[t]}=\frac{\beta ^{\lbrack t]}-\frac{\sum_{\tau =t}^{\tau
=L}\beta ^{\lbrack \tau ]}}{L-t+1}}{2\left( \varepsilon +\gamma \right) }+%
\frac{\left( 1-\frac{\sum_{\tau =t}^{\tau =L}\delta ^{\tau -t}}{L-t+1}%
\right) \mathbb{E}\left[ \left. \eta ^{\lbrack t]}\right\vert X^{[t]}\right]
}{2\left( \varepsilon +\gamma \right) }-\frac{\sum_{\tau =1}^{\tau
=t-1}d^{[\tau ]}}{L-t+1}.
\end{equation}%
$\square $

\textbf{Proof of Proposition \ref{prop_two period}}. The payoff can be
simplified by plugging $d_{i}^{[2]}=-d_{i}^{[1]}.$ To derive the equilibrium
storage quantities, we set $\frac{\partial \mathbb{E}[\pi _{i}|x_{i}^{[1]}]}{%
\partial d_{i}^{[1]}}=0$:

\begin{equation}
\left[
\begin{array}{c}
\beta ^{\lbrack 1]}-\beta ^{\lbrack 2]}-\left( \gamma ^{\lbrack 1]}+\gamma
^{\lbrack 2]}\right) \sum_{j\neq i}\mathbb{E}[d_{j}^{[1]}|x_{i}^{[1]}] \\
+\mathbb{E}[\eta ^{\lbrack 1]}|x_{i}^{[1]}]-\mathbb{E}[\eta ^{\lbrack
2]}|x_{i}^{[1]}] \\
-2\left( \gamma ^{\lbrack 1]}+\gamma ^{\lbrack 2]}+\varepsilon ^{\lbrack
1]}+\varepsilon ^{\lbrack 2]}\right) \cdot d_{i}^{[1]}%
\end{array}%
\right] =0.
\end{equation}%
Notice that

\begin{equation*}
\mathbb{E}[d_{j}^{[1]}|x_{i}^{[1]}]=A+Bx_{0}^{[1]}+C\mathbb{E}%
[x_{j}^{[1]}|,x_{i}^{[1]}],
\end{equation*}

\begin{eqnarray*}
\mathbb{E}[x_{j}^{[1]}|x_{i}^{[1]}] &=&\mathbb{E}[\eta ^{\lbrack 1]}+\xi
_{j}^{[1]}|x_{i}^{[1]}] \\
&=&\mathbb{E}[\eta ^{\lbrack 1]}|x_{i}^{[1]}],
\end{eqnarray*}

\begin{eqnarray*}
\mathbb{E}[\eta ^{\lbrack 2]}|x_{i}^{[1]}] &=&\mathbb{E}[\delta \eta
^{\lbrack 1]}+\epsilon _{1}|x_{i}^{[1]}] \\
&=&\delta \mathbb{E}[\eta ^{\lbrack 1]}|x_{i}^{[1]}],
\end{eqnarray*}%
\begin{equation*}
\mathbb{E}[\eta ^{\lbrack 1]}|x_{i}^{[1]}]=\frac{\rho }{\alpha +\rho }%
x_{i}^{[1]}.
\end{equation*}%
By matching the coefficients with respect to $x_{i}^{[1]}$, we have

\begin{equation*}
A=\frac{\beta ^{\lbrack 1]}-\beta ^{\lbrack 2]}}{2\left( \varepsilon
^{\lbrack 1]}+\varepsilon ^{\lbrack 2]}\right) +(n+1)\left( \gamma ^{\lbrack
1]}+\gamma ^{\lbrack 2]}\right) },
\end{equation*}

\begin{equation*}
C=\frac{(1-\delta )\rho }{\left[
\begin{array}{c}
(n-1)\left( \gamma ^{\lbrack 1]}+\gamma ^{\lbrack 2]}\right) \rho + \\
2\left( \varepsilon ^{\lbrack 1]}+\varepsilon ^{\lbrack 2]}+\gamma ^{\lbrack
1]}+\gamma ^{\lbrack 2]}\right) (\alpha +\rho )%
\end{array}%
\right] }.
\end{equation*}
$\square$

\textbf{Proof of Proposition \ref{prop:twoperiodpayoff}}. The payoff can be
calculated through
\begin{eqnarray}
\mathbb{E}[\pi _{i}] &=&\mathbb{E}\left[ \mathbb{E}[\pi _{i}|x_{i}^{[1]}]%
\right]  \notag \\
&=&\left( \varepsilon ^{\lbrack 1]}+\varepsilon ^{\lbrack 2]}+\gamma
^{\lbrack 1]}+\gamma ^{\lbrack 2]}\right)  \notag \\
&&\cdot \left( A^{2}+C^{2}\mathbb{E}[x_{i}^{[1]}]^{2}\right) .
\end{eqnarray}%
Notice that storage $i$'s payoff $\mathbb{E}[\pi _{i}]=\sum_{t\in T}\left(
\varepsilon ^{\lbrack t]}+\gamma ^{\lbrack t]}\right) A^{2}$ when there is
no information available. The additional payoff proportional to $C^{2}%
\mathbb{E}[x_{i}^{[1]}]^{2}$ corresponds to the economic value of the
private forecast.

\begin{eqnarray}
\lim_{n\rightarrow \infty }\mathbb{E}[\pi _{i}] &=&\lim_{n\rightarrow \infty
}\sum_{t=1,2}\left( \varepsilon ^{\lbrack t]}+\gamma ^{\lbrack t]}\right) %
\left[
\begin{array}{c}
\frac{\left( \beta ^{\lbrack 1]}-\beta ^{\lbrack 2]}\right) ^{2}}{%
(n+1)^{2}\left( \gamma ^{\lbrack 1]}+\gamma ^{\lbrack 2]}\right) 2} \\
+\frac{(1-\delta )^{2}}{(n-1)^{2}\left( \gamma ^{\lbrack 1]}+\gamma
^{\lbrack 2]}\right) ^{2}\rho }%
\end{array}%
\right]  \notag \\
&&\lim_{n\rightarrow \infty }\sum_{t=1,2}\left( \varepsilon ^{\lbrack
t]}+\gamma ^{\lbrack t]}\right) \left[
\begin{array}{c}
\left( \frac{\beta ^{\lbrack 1]}-\beta ^{\lbrack 2]}}{\gamma ^{\lbrack
1]}+\gamma ^{\lbrack 2]}}\right) ^{2}+ \\
\left( \frac{1-\delta }{\gamma ^{\lbrack 1]}+\gamma ^{\lbrack 2]}}\right)
^{2}\rho ^{-1}%
\end{array}%
\right] \cdot n^{-2}.
\end{eqnarray}%
$\square $

\textbf{Proof of Proposition \ref{prop_public}}. To derive the equilibrium
storage quantities, we set $\frac{\partial \mathbb{E}[\pi _{i}|x_{0}^{[1]}]}{%
\partial d_{i}^{[1]}}=0$:

\begin{equation}
\left[
\begin{array}{c}
\beta ^{\lbrack 1]}-\beta ^{\lbrack 2]}-\left( \gamma ^{\lbrack 1]}+\gamma
^{\lbrack 2]}\right) \sum_{j\neq i}\mathbb{E}[d_{j}^{[1]}|x_{0}^{[1]}] \\
+\mathbb{E}[\eta ^{\lbrack 1]}|x_{0}^{[1]}]-\mathbb{E}[\eta ^{\lbrack
2]}|x_{0}^{[1]}] \\
-2\left( \gamma ^{\lbrack 1]}+\gamma ^{\lbrack 2]}+\varepsilon ^{\lbrack
1]}+\varepsilon ^{\lbrack 2]}\right) \cdot d_{i}^{[1]}%
\end{array}%
\right] =0.
\end{equation}%
Notice that $\mathbb{E}[d_{j}^{[1]}|x_{0}^{[1]}]=A+Bx_{0}^{[1]},$ since $%
x_{0}^{[1]}$ is common knowledge. $\mathbb{E}[\eta ^{\lbrack
1]}|x_{0}^{[1]}]=\frac{\sigma }{\alpha +\sigma }x_{0}^{[1]}$, and

\begin{eqnarray*}
\mathbb{E}[\eta ^{\lbrack 2]}|x_{0}^{[1]}] &=&\mathbb{E}[\delta \eta
^{\lbrack 1]}+\epsilon _{1}|x_{0}^{[1]}] \\
&=&\delta \mathbb{E}[\eta ^{\lbrack 1]}|x_{0}^{[1]}].
\end{eqnarray*}%
By matching the coefficients with respect to $x_{i}^{[1]}$, we have

\begin{equation*}
A=\frac{\beta ^{\lbrack 1]}-\beta ^{\lbrack 2]}}{2\left( \varepsilon
^{\lbrack 1]}+\varepsilon ^{\lbrack 2]}\right) +(n+1)\left( \gamma ^{\lbrack
1]}+\gamma ^{\lbrack 2]}\right) },
\end{equation*}

\begin{equation*}
B=\frac{(1-\delta )\sigma /(\alpha +\sigma )}{\left[
\begin{array}{c}
2\left( \varepsilon ^{\lbrack 1]}+\varepsilon ^{\lbrack 2]}+\gamma ^{\lbrack
1]}+\gamma ^{\lbrack 2]}\right) \\
+\left( \gamma ^{\lbrack 1]}+\gamma ^{\lbrack 2]}\right) (n-1)%
\end{array}%
\right] }.
\end{equation*}%
The corresponding payoff can be calculated as follows:

\begin{equation*}
\mathbb{E}[\pi _{i}]=\sum_{t=1,2}\left( \varepsilon ^{\lbrack t]}+\gamma
^{\lbrack t]}\right)
\end{equation*}%
\begin{equation}
\cdot \left\{
\begin{array}{c}
\left[ \frac{\beta ^{\lbrack 1]}-\beta ^{\lbrack 2]}}{\sum_{t=1,2}2\left(
\varepsilon ^{\lbrack t]}+\gamma ^{\lbrack t]}\right) +(n-1)\left( \gamma
^{\lbrack 1]}+\gamma ^{\lbrack 2]}\right) }\right] ^{2} \\
+\left[ \frac{(1-\delta )}{\sum_{t=1,2}2\left( \varepsilon ^{\lbrack
t]}+\gamma ^{\lbrack t]}\right) +\left( \gamma ^{\lbrack 1]}+\gamma
^{\lbrack 2]}\right) (n-1)}\right] ^{2}\cdot \frac{\sigma }{\left( \alpha
+\sigma \right) ^{2}}%
\end{array}%
\right\} .
\end{equation}%
\begin{equation*}
\lim_{n\rightarrow \infty }\mathbb{E}[\pi _{i}]\rightarrow
\sum_{t=1,2}\left( \varepsilon ^{\lbrack t]}+\gamma ^{\lbrack t]}\right) %
\left[
\begin{array}{c}
\left( \frac{\beta ^{\lbrack 1]}-\beta ^{\lbrack 2]}}{\gamma ^{\lbrack
1]}+\gamma ^{\lbrack 2]}}\right) ^{2} \\
+\left( \frac{1-\delta }{\gamma ^{\lbrack 1]}+\gamma ^{\lbrack 2]}}\right)
^{2}\cdot \frac{\sigma }{\left( \alpha +\sigma \right) ^{2}}%
\end{array}%
\right] \cdot n^{-2}.
\end{equation*}%
$\square $

\textbf{Proof of Proposition \ref{prop_target}}. To solve for equilibrium
storage quantities, we set $\frac{\partial \mathbb{E}[\pi _{i}|x_{0}^{[1]}]}{%
\partial d_{i}^{[1]}}=0$ for $\forall i$ $\in J$ and $\frac{\partial \pi _{i}%
}{\partial d_{i}^{[1]}}=0$ for $\forall i$ $\in I-J$, separately. For $%
\forall i\in J,$

\begin{equation}
\begin{array}{c}
\beta ^{\lbrack 1]}-\beta ^{\lbrack 2]}-\left( \gamma ^{\lbrack 1]}+\gamma
^{\lbrack 2]}\right) \left[
\begin{array}{c}
(n-m)C+ \\
(m-1)\left( A+Bx_{0}^{[1]}\right)%
\end{array}%
\right] \\
+\frac{(1-\delta )\sigma }{\alpha +\sigma }x_{0}^{[1]}-2\left(
\begin{array}{c}
\gamma ^{\lbrack 1]}+\gamma ^{\lbrack 2]} \\
+\varepsilon ^{\lbrack 1]}+\varepsilon ^{\lbrack 2]}%
\end{array}%
\right) \cdot \left( A+Bx_{0}^{[1]}\right)%
\end{array}%
=0,
\end{equation}%
whereas for $\forall i\in I-J,$

\begin{equation}
\begin{array}{c}
\beta ^{\lbrack 1]}-\beta ^{\lbrack 2]}-\left( \gamma ^{\lbrack 1]}+\gamma
^{\lbrack 2]}\right) \left[
\begin{array}{c}
mA+ \\
(n-m-1)C%
\end{array}%
\right] \\
-2\left(
\begin{array}{c}
\gamma ^{\lbrack 1]}+\gamma ^{\lbrack 2]} \\
+\varepsilon ^{\lbrack 1]}+\varepsilon ^{\lbrack 2]}%
\end{array}%
\right) \cdot C%
\end{array}%
=0.
\end{equation}%
Matching coefficients with respect to $x_{0}^{[1]}$, we can obtain $A,B$ and
$C$ following a similar procedure as before. We measure the economic
efficiency by aggregate payoff:

\begin{equation}
\sum_{i\in I}\mathbb{E}[\pi _{i}]=\left( \varepsilon ^{\lbrack
1]}+\varepsilon ^{\lbrack 2]}+\gamma ^{\lbrack 1]}+\gamma ^{\lbrack
2]}\right)  \notag
\end{equation}%
\begin{equation}
\cdot \left( nA^{2}+\frac{(1-\delta )^{2}m}{\left[
\begin{array}{c}
2\left( \gamma ^{\lbrack 1]}+\gamma ^{\lbrack 2]}+\varepsilon ^{\lbrack
1]}+\varepsilon ^{\lbrack 2]}\right) \\
+(m-1)\left( \gamma ^{\lbrack 1]}+\gamma ^{\lbrack 2]}\right)%
\end{array}%
\right] ^{2}}\frac{\sigma }{\left( \alpha +\sigma \right) ^{2}}\right) .
\end{equation}%
It can be checked that the storages' aggregate payoff $\sum_{i\in I}\mathbb{E%
}[\pi _{i}]$ is maximized when%
\begin{equation}
m=1+\frac{2\left( \varepsilon ^{\lbrack 1]}+\varepsilon ^{\lbrack 2]}\right)
}{\gamma ^{\lbrack 1]}+\gamma ^{\lbrack 2]}}.
\end{equation}%
$\square $

\textbf{Proof of Proposition \ref{prop_multiple period}}. The payoff from
storage $i$ under a Lagrangian relaxation can be expressed as%
\begin{equation*}
L_{i}\left( d_{i}^{[t]},t\in T\right) =\sum_{t\in T}l_{i}^{[t]}\left(
d_{i}^{[t]}\right)
\end{equation*}

\begin{eqnarray}
&=&\sum_{t\in T}\mathbb{E}\left[ \left.
\begin{array}{c}
P^{[t]}(D^{[t]})\cdot d_{i}^{[t]} \\
-\varepsilon ^{\lbrack t]}\cdot \left( d_{i}^{[t]}\right) ^{2}%
\end{array}%
\right\vert X_{i}^{[t]}\right] -\lambda ^{\lbrack t]}\mathbb{E}\left[ \left.
\sum_{t\in T}d_{i}^{[t]}\right\vert X_{i}^{[t]}\right]   \notag \\
&=&\sum_{t\in T}\mathbb{E}\left[ \left.
\begin{array}{c}
P^{[t]}(D^{[t]})\cdot d_{i}^{[t]} \\
-\varepsilon ^{\lbrack t]}\cdot \left( d_{i}^{[t]}\right) ^{2}-\lambda
^{\lbrack t]}d_{i}^{[t]}%
\end{array}%
\right\vert X_{i}^{[t]}\right] ,
\end{eqnarray}%
wherein $\lambda ^{\lbrack t]}$ is a Lagrangian multiplier to relax the
constraint that $\mathbb{E}\left[ \left. \sum_{t\in T}d_{i}^{[t]}\right\vert
X_{i}^{[t]}\right] =0$. To solve for equilibrium storage quantities, we set $%
\frac{\partial l_{i}^{[t]}\left( d_{i}^{[t]}\right) }{\partial d_{i}^{[t]}}=0
$:

\begin{equation}
\left[
\begin{array}{c}
\beta ^{\lbrack t]}-\gamma ^{\lbrack t]}\sum_{j\neq i}\mathbb{E}%
[d_{j}^{[t]}|X_{i}^{[t]}] \\
+\mathbb{E}[\eta ^{\lbrack t]}|X_{i}^{[t]}]-2\left( \gamma ^{\lbrack
t]}+\varepsilon ^{\lbrack t]}\right) \cdot d_{i}^{[1]}-\lambda ^{\lbrack t]}%
\end{array}%
\right] =0.
\end{equation}%
By $\mathbb{E}\left[ \sum_{t\in T}\left( d_{i}^{[t]}\right) ^{\ast }\right]
=0$, we use $\lambda $ as a static approximate solution instead of $\lambda
^{\lbrack t]}$. $\square $

\textbf{Proof of Proposition \ref{prop_Hetero}}. Again, the first-order
condition for payoff-maximization requires

\begin{equation}
\left[
\begin{array}{c}
\beta ^{\lbrack 1]}-\beta ^{\lbrack 2]}-\left( \gamma ^{\lbrack 1]}+\gamma
^{\lbrack 2]}\right) \sum_{j\neq i}\mathbb{E}%
[d_{j}^{[1]}|x_{0}^{[1]},x_{i}^{[1]}] \\
+\mathbb{E}[\eta ^{\lbrack 1]}|x_{0}^{[1]},x_{i}^{[1]}]-\mathbb{E}[\eta
^{\lbrack 2]}|x_{0}^{[1]},x_{i}^{[1]}] \\
-2\left( \gamma ^{\lbrack 1]}+\gamma ^{\lbrack 2]}+\varepsilon
_{i}^{[1]}+\varepsilon _{i}^{[2]}\right) \cdot d_{i}^{[1]}%
\end{array}%
\right] =0.
\end{equation}%
The unknown coefficients in the equilibrium buying or selling quantities are
summarized as follows.
\begin{eqnarray}
A_{i} &=&\frac{\beta ^{\lbrack 1]}-\beta ^{\lbrack 2]}}{2\left( \varepsilon
_{i}^{[1]}+\varepsilon _{i}^{[2]}\right) +\left( \gamma ^{\lbrack 1]}+\gamma
^{\lbrack 2]}\right) }  \notag \\
&&\cdot \left[ \sum_{i\in I}\frac{\left( \gamma ^{\lbrack 1]}+\gamma
^{\lbrack 2]}\right) }{2\left( \varepsilon _{i}^{[1]}+\varepsilon
_{i}^{[2]}\right) +\left( \gamma ^{\lbrack 1]}+\gamma ^{\lbrack 2]}\right) }%
+1\right] ^{-1}.
\end{eqnarray}

\begin{eqnarray}
C_{i} &=&-\frac{\left( \gamma ^{\lbrack 1]}+\gamma ^{\lbrack 2]}\right)
\frac{(1-\delta )\rho _{i}}{\alpha +\sigma +\rho _{i}}}{\left[ 2\left(
\begin{array}{c}
\gamma ^{\lbrack 1]}+\gamma ^{\lbrack 2]} \\
+\varepsilon _{i}^{[1]}+\varepsilon _{i}^{[2]}%
\end{array}%
\right) -\left( \gamma ^{\lbrack 1]}+\gamma ^{\lbrack 2]}\right) \frac{\rho
_{i}}{\alpha +\sigma +\rho _{i}}\right] }  \notag \\
&&\cdot \frac{\sum_{i\in I}\frac{\frac{\rho _{i}}{\alpha +\sigma +\rho _{i}}%
}{\left[ 2\left( \gamma ^{\lbrack 1]}+\gamma ^{\lbrack 2]}+\varepsilon
_{i}^{[1]}+\varepsilon _{i}^{[2]}\right) -\left( \gamma ^{\lbrack 1]}+\gamma
^{\lbrack 2]}\right) \frac{\rho _{i}}{\alpha +\sigma +\rho _{i}}\right] }}{%
1+\sum_{i\in I}\frac{\left( \gamma ^{\lbrack 1]}+\gamma ^{\lbrack 2]}\right)
\frac{\rho _{i}}{\alpha +\sigma +\rho _{i}}}{\left[ 2\left( \gamma ^{\lbrack
1]}+\gamma ^{\lbrack 2]}+\varepsilon _{i}^{[1]}+\varepsilon
_{i}^{[2]}\right) -\left( \gamma ^{\lbrack 1]}+\gamma ^{\lbrack 2]}\right)
\frac{\rho _{i}}{\alpha +\sigma +\rho _{i}}\right] }}  \notag \\
&&+\frac{\frac{(1-\delta )\rho _{i}}{\alpha +\sigma +\rho _{i}}}{\left[
2\left(
\begin{array}{c}
\gamma ^{\lbrack 1]}+\gamma ^{\lbrack 2]} \\
+\varepsilon _{i}^{[1]}+\varepsilon _{i}^{[2]}%
\end{array}%
\right) -\left( \gamma ^{\lbrack 1]}+\gamma ^{\lbrack 2]}\right) \frac{\rho
_{i}}{\alpha +\sigma +\rho _{i}}\right] }.
\end{eqnarray}

\begin{eqnarray}
B_{i} &=&-\frac{\left( \gamma ^{\lbrack 1]}+\gamma ^{\lbrack 2]}\right)
\left( \sum_{i\in I}B_{i}+\frac{\sigma }{\alpha +\sigma +\rho _{i}}%
\sum_{j\neq i}C_{j}\right) }{\left[ 2\left( \varepsilon
_{i}^{[1]}+\varepsilon _{i}^{[2]}\right) +\left( \gamma ^{\lbrack 1]}+\gamma
^{\lbrack 2]}\right) \right] }  \notag \\
&&+\frac{\frac{(1-\delta )\sigma }{\alpha +\sigma +\rho _{i}}}{\left[
2\left( \varepsilon _{i}^{[1]}+\varepsilon _{i}^{[2]}\right) +\left( \gamma
^{\lbrack 1]}+\gamma ^{\lbrack 2]}\right) \right] },
\end{eqnarray}%
wherein
\begin{eqnarray}
\sum_{i\in I}B_{i} &=&-\frac{\sum_{i\in I}\frac{\left( \gamma
_{i}^{[1]}+\gamma _{i}^{[2]}\right) \frac{\sigma }{\alpha +\sigma +\rho _{i}}%
\sum_{j\neq i}C_{j}}{2\left( \varepsilon _{i}^{[1]}+\varepsilon
_{i}^{[2]}\right) +\left( \gamma _{i}^{[1]}+\gamma _{i}^{[2]}\right) }}{%
\left[ 1+\sum_{i\in I}\frac{\left( \gamma _{i}^{[1]}+\gamma
_{i}^{[2]}\right) }{2\left( \varepsilon _{i}^{[1]}+\varepsilon
_{i}^{[2]}\right) +\left( \gamma _{i}^{[1]}+\gamma _{i}^{[2]}\right) }\right]
}  \notag \\
&&+\frac{\sum_{i\in I}\frac{\frac{(1-\delta )\sigma }{\alpha +\sigma +\rho
_{i}}}{2\left( \varepsilon _{i}^{[1]}+\varepsilon _{i}^{[2]}\right) +\left(
\gamma ^{\lbrack 1]}+\gamma ^{\lbrack 2]}\right) }}{\left[ 1+\sum_{i\in I}%
\frac{\left( \gamma ^{\lbrack 1]}+\gamma ^{\lbrack 2]}\right) }{2\left(
\varepsilon _{i}^{[1]}+\varepsilon _{i}^{[2]}\right) +\left( \gamma
^{\lbrack 1]}+\gamma ^{\lbrack 2]}\right) }\right] }.
\end{eqnarray}
$\square$

\end{document}